\documentclass[11pt]{article}
\usepackage{amsfonts}
\usepackage{latexsym,amsmath}
\usepackage{amssymb,array}
\usepackage{epsfig}
\usepackage{color}
\makeatletter \oddsidemargin 0in \evensidemargin 0in \textwidth
16cm
\newtheorem{t1}{Theorem}[section]

\newtheorem{l1}{Lemma}[section]
\newtheorem{c1}{Corollary}[section]

\newtheorem{r1}{Remark}[section]
\newtheorem{co}{Counterexample}[section]
\newtheorem{ex}{Example}[section]
\newcommand{\bea}{\begin{eqnarray}}
\newcommand{\eea}{\end{eqnarray}}

\begin{document}
\title{Some Ageing Properties of Dynamic Additive Mean Residual Life Model\footnote{This article has been published in ``Rashi (2017), Vol. 2, Issue 1, pp. 26-33''}}
\author{Suchismita Das\\303/304, Shri Raghunarayan CHS\\
	MIDC Road\\
	Thakurli(E), Mumbai
\and Asok K. Nanda\footnote{corresponding author: asok@iiserkol.ac.in; asok.k.nanda@gmail.com}\\
Department of Mathematics and Statistics\\
IISER Kolkata}
\date{March 2017}
\maketitle
\begin{abstract}
Although proportional hazard rate model is a very popular model to analyze failure time data,
sometimes it becomes important to study the additive hazard rate model. Again, sometimes the
concept of the hazard rate function is abstract, in comparison to the concept of mean residual
life function. A new model called `dynamic additive mean
residual life model' where the covariates are time-dependent has been defined in the literature. Here we study the closure properties of
the model for different positive and negative ageing classes under certain condition(s). Quite a few examples are presented to illustrate
different properties of the model.
\end{abstract}
{\bf Key Words and Phrases:} Additive mean residual life function, ageing property.\\
{\bf AMS 2010 Classifications:} Primary 62N05; Secondary 90B25.
\section{Introduction}
\hspace{0.5in}In the literature, a large number of papers deal with modeling and analyzing data
on the time until occurrence of an event. These events sometimes represent the time
to failure of a system of components or a living organism. Cox's \cite{c} proportional
hazard rate model has been used to model failure time data. This model is very popular for
analysis of right-censored data, and has been used for estimating the risk of failure
associated with vector of covariates.

\hspace{0.5in}Although multiplicative hazards model is mostly studied in literature, it is
important to study the additive hazard rate model as well. This is because from the public health point
of view it is very important to study the risk difference than
the risk ratio in describing the association between the risk factor and the occurrence of a
disease, see, for instance, Breslow and Day \cite{bd,brd}. They have shown that additive hazard
rate model fits certain type of data better than
the proportional hazard model. Lin and Ying \cite{ly} have shown that the additive hazards model
provides a simple structure for studying frailty models and interval-censored data, which are very
difficult to deal with under the proportional hazards model. A family of proportional and of additive
hazards models for the analysis of grouped survival data has been considered by Tibshirani and Ciampi
\cite{tc}. Mckeague and Utikal \cite{mu} have developed goodness-of-fit tests for Cox's proportional hazards
model and Aalen's additive risk model, and each model has been compared on an equal footing with the
best fitted fully nonparametric model. Using generalized linear models (GLM), Hakulinen and Tenkanen \cite{ht} have shown
how a proportional hazards regression model may be adapted to the relative survival rates.
Bad$\acute{i}$a, Berrade and Campos \cite{bbc} have studied
some ageing characteristics of additive and of proportional hazard mixing models. They
have also studied the effect of mixing on stochastic ordering. Martinussen and Scheike \cite{ms} have
compared the full Aalen additive hazards model and the change-point model, and discussed
how to estimate the parameters of the change-point model. The analytical properties of additive
hazards model have been studied by Nair and Sankaran \cite{ns}. They have compared the ageing
properties of the baseline random variable and the induced random variable. Li and Ling \cite{ll}
have discussed the ageing and the dependence properties in the additive hazard mixing model. Some
stochastic comparisons have also been studied in this paper. Bin \cite{b} has discussed regression
analysis of failure time under the additive hazards model, when the regression coefficients
are time-varying.

\hspace{0.5in}Sometimes the concept of hazard function becomes abstract, in comparison to
the concept of mean residual life (MRL) function. The hazard rate
is the instantaneous failure rate at any point of time,
whereas the MRL summarizes the entire residual life. The MRL function has more
intuitive appeal for modeling and analysis of failure data than the concept of hazard rate
function. With this in mind, the additive mean residual life (AMRL) model has been developed.
The AMRL model specifies that the MRL function associated with the covariates is the sum of the baseline
MRL function and the constant representing the function of covariates. For example, a new drug prescribed to a patient
may work well in the beginning of the treatment period, but after a certain time the effect of the drug may decrease.
Then it is very important to know when and how fast the drug becomes ineffective, see, for instance, Bin \cite{b}. AMRL model helps us design further studies
to explore the treatment strategies for patients. In some practical situations, the covariates may not be constant over the
whole time interval $[0,\infty)$, but they may vary over different time intervals. Using this idea Das and Nanda \cite{dn} have developed a new model called
dynamic additive mean residual life (DAMRL) model. They have studied the closure of this model under different stochastic orders.

\hspace{0.5in}In Section~$2$ we give a brief description of the DAMRL model, Section~$3$ provides
some ageing classes under the DAMRL model with some illustrative examples, and Section $4$ concludes the manuscript.

\hspace{0.5in}Throughout the paper, {\it increasing} and {\it decreasing}
properties of a function are not used in strict sense. For any
twice differentiable function $g(t)$, we write $g'(t)$ and
$g''(t)$ to denote the first and the second derivatives of
$g(t)$ with respect to $t$, respectively. We denote by $a\stackrel{\rm
sign}=b$ to mean that $a$ and $b$ have the same sign, whereas $\alpha\stackrel{\rm
def}=\beta$ tells that $\alpha$ is defined by $\beta$.

\section{Dynamic Additive Mean Residual Life Model}
\hspace{0.5in}Let $X$ be a nonnegative random variable with finite mean, survival function $\bar{F}_{X}(\cdot)$ and hazard rate function $r_{X}(\cdot)$. Then the MRL is given by
\begin{eqnarray*}
m_{X}(t)=E\left[X-t|X>t\right]=\frac{\int_{t}^{\infty}\bar{F}_{X}(u)du}{\bar{F}_{X}(t)},
\end{eqnarray*}
for $t\geqslant 0$ such that $\bar{F}_{X}(t)>0$.

\hspace{0.5in}Let $X^*$ be a nonnegative absolutely continuous random variables
having MRL function $m_{X^*}(\cdot)$. Then the AMRL model is given by
\begin{eqnarray}\label{e4}
m_{X^*}(t)&=&c+m_{X}(t),
\end{eqnarray}
where $c$ is a function of covariates, independent of $t$ such that $c+m_{X}(t)$ is nonnegative for all $t$. Das and Nanda \cite{dn} have shown
that Additive Hazard Rate model and AMRL model in general do not imply one another. For some more work in this direction, one may refer to Yin and Cai \cite{yc}, Sun and Zhang \cite{sz} and the references there in.
If the covariates are time-dependent, the AMRL model reduces to dynamic AMRL model given by
\begin{eqnarray}\label{e1}
m_{X^*}(t)&=&c(t)+m_{X}(t).
\end{eqnarray}
Before going to discuss some results we give a lemma without proof, which gives some conditions on c(t) [cf. Das and Nanda \cite{dn}].
\begin{l1}\label{l1}
If two nonnegative random variables $X$ and $X^*$ satisfy (\ref{e1}), then the following conditions must be satisfied:
\begin{enumerate}
\item[(i)]$~0\leqslant c(t)+m_{X}(t)< \infty,~~$ for all $~t\geqslant 0$;
\item[(ii)]$~c(t)$ is a continuous function of $~t\geqslant 0$;
\item[(iii)]$~t+c(t)+m_{X}(t)$ is increasing in $~t\geqslant 0$;
\item[(iv)]if there does not exist any $t_{0}$ with $m_{X}(t_{0})=0$, then
$\displaystyle\int_{0}^{\infty}\frac{dt}{c(t)+m_{X}(t)}=\infty$. \hfill$\Box$
\end{enumerate}
\end{l1}
\begin{r1}
It is easy to see that if $c(t)$ is increasing in $t\geqslant 0$, then the condition ($iii$) of Lemma~\ref{l1} trivially holds, whereas when $c(t)\geqslant 0$, for all $t$, condition ($i$) of Lemma~\ref{l1} automatically holds.\hfill$\Box$
\end{r1}
Again, it can be noted that
\begin{eqnarray}\label{e2}
r_{X^*}(t)&=&\frac{1+m_{X^*}'(t)}{m_{X^*}(t)}\nonumber\\
&=&r_{X}(t)\frac{m_{X}(t)}{c(t)+m_{X}(t)}+\frac{c'(t)}{c(t)+m_{X}(t)}.
\end{eqnarray}
\section{Properties of Some Ageing Classes}
Ageing is an inherent property of a unit that may be a system of components or a living organism. It is characterized by various quantities, viz., hazard rate,
mean residual life etc. By ageing we generally mean the adverse effect of age on the random residual lifetime of a unit. To be specific, by ageing we generally mean positive ageing, which means that an older system has a shorter remaining lifetime, in some statistical sense, than a newer one. In this section we have studied the closure of different ageing properties such as IFR (increasing in failure rate), IFRA (increasing failure rate in average), NBU (new better than used), NBUFR (new better than used in failure rate), NBAFR (new better than used in failure rate average) and their duals with time-dependent covariate(s).

\hspace{0.5in}In this section we study the condition(s) under which $X$ and $X^*$ satisfying (\ref{e1}) share some ageing properties. The following
theorem shows that under the model (\ref{e1}), the IFR (resp.
DFR) property of $X$ is transmitted to the random variable
$X^{*}$ under certain condition on c(t). Keep in mind that a random variable $X$ with failure rate function $r_{X}(\cdot)$ is said to be IFR  (resp. DFR) if $r_{X}(t)$ is
increasing (resp. decreasing) in $t$.

\hspace{0.5in}The following theorem states the conditions for $X^{*}$ to be IFR. Here we need $c(t)$ to be positive.
\begin{t1}\label{th1}
If the random variable $X$ is IFR, then the random variable $X^*$ satisfying (\ref{e1}) is IFR provided, for all $t\geqslant 0$,
\begin{enumerate}
\item[(i)] $~\frac{m_{X}(t)}{c(t)}$ is increasing in $t$;
\item[(ii)]$~c(t)$ is logconvex.
\end{enumerate}
\end{t1}
$\mathbf{Proof}$: Differentiating (\ref{e2}) with respect to $t$, we get
\begin{eqnarray}\label{e6}
r_{X^*}'(t)&=&r_{X}'(t)\frac{m_{X}(t)}{c(t)+m_{X}(t)}+r_{X}(t)\frac{m_{X}'(t)c(t)-m_{X}(t)c'(t)}{\left(c(t)+m_{X}(t)\right)^2}
+\frac{c(t)c^{''}(t)-c'^2(t)}{\left(c(t)+m_{X}(t)\right)^2}\nonumber\\
&&\nonumber\\
&&+\frac{m_{X}(t)c^{''}(t)-c'(t)m_{X}'(t)}{\left(c(t)+m_{X}(t)\right)^2}.
\end{eqnarray}
Now, the first term in the above expression is nonnegative because $X$ is IFR, the second term is nonnegative due to ($i$) in the hypothesis, the third term is nonnegative because of ($ii$), and finally the last term is nonnegative if $\frac{c'(t)}{m_{X}(t)}$ is increasing in $t$, which is always true. This is because, by ($ii$), $c(t)$ is convex, which implies that $c'(t)$ is increasing in $t$, and since $X$ is IFR, it is also DMRL, so that $m_{X}(t)$ is decreasing in $t$.
Hence, $X^*$ is IFR.~\hfill$\Box$

\hspace{0.5in}We now present an application of Theorem \ref{th1}.
\begin{ex}\label{ex13}
Let $X$ follow standard exponential distribution. Take $c(t)=\exp(-t)$, for all $t\geqslant~0$.
Clearly, $c(t)$
satisfies all the conditions of Lemma~\ref{l1}. Again, conditions ($i$) and ($ii$) of Theorem~\ref{th1} are also satisfied.
Hence, by Theorem~\ref{th1}, $X^*$ is IFR.~\hfill$\Box$
\end{ex}
\hspace{0.5in}The following counterexample shows that condition
($i$) of Theorem \ref{th1} is a sufficient condition but not necessary.
\begin{co}\label{co4}
Let $X$ be a random variable having mean residual life $m_{X}(t)=1/(2+t)$, $t\geqslant 0$. Take $c(t)=1/(3+t)$, for all $t\geqslant 0$.
Then $c(t)$ satisfies all the conditions of Lemma \ref{l1}. Again, condition ($ii$) of Theorem \ref{th1} is satisfied while condition ($i$) is not. Now, for all $t\geqslant 0$,
\begin{eqnarray*}
r_{X^*}(t)&=&\frac{((2+t)^2-1)(3+t)}{(2+t)(5+2t)}-\frac{(2+t)}{(3+t)(5+2t)},\\
\end{eqnarray*}
which can be shown to be increasing in $t$. Thus, $X^{*}$ is IFR. Hence, condition
($i$) of Theorem~\ref{th1} is a sufficient condition but not necessary.~\hfill$\Box$
\end{co}
\begin{r1}\label{r1}
By taking $X$ to be a standard exponential distribution, and
\begin{eqnarray}
c(t)=
  \begin{cases}
      1/(2+t)^2, ~~ ~~~~~~0\leqslant t\leqslant 1, &
     \\
      1/(3(2+t^2)),~~~~~~~~~~t \geqslant 1,\nonumber
  \end{cases}
\end{eqnarray}
one can show that condition ($ii$) of Theorem \ref{th1} is a sufficient condition but not necessary.~\hfill$\Box$
\end{r1}
\hspace{0.5in}The following theorem states the conditions for $X^{*}$ to be DFR. Here we need $c(t)$ to be positive. The proof being similar to that of Theorem \ref{th1} is omitted.
\begin{t1}\label{th2}
If the random variable $X$ is DFR, then the random variable $X^*$ satisfying (\ref{e1}) is DFR provided,
\begin{enumerate}
\item[(i)] $~\frac{m_{X}(t)}{c(t)}$ is decreasing in $t$;
\item[(ii)]$~c(t)$ is increasing and concave;
\end{enumerate}
for all $t\geqslant 0$. \hfill$\Box$
\end{t1}

\hspace{0.5in}Here we present an application of Theorem \ref{th2}.
\begin{ex}\label{ex1}
Let $X$ be a random variable having the failure rate
\begin{eqnarray*}
r_{X}(t)=
  \begin{cases}
     2/(1+t), ~~ ~0\leqslant t\leqslant 1, &
     \\
     1, ~~~~ ~~ ~~~~~~~~~~~t \geqslant 1.
  \end{cases}
 \end{eqnarray*}
Take $c(t)=t/(1+t)$, for all $t\geqslant 0$.
Clearly, $c(t)$ satisfies all the conditions of Lemma~\ref{l1}. One can verify that $\frac{m_{X}(t)}{c(t)}$ is decreasing in $t\geqslant 0$ and also $c(t)$ is concave. Thus, by Theorem~\ref{th2}, $X^*$ is DFR.~\hfill$\Box$
\end{ex}
\begin{r1}\label{r2}
In Example \ref{ex1} if we take $c(t)=\frac{1+t}{2+t}$, for all $t\geqslant 0$, then it can be seen that condition ($i$)
of Theorem~\ref{th2} is a sufficient condition but not necessary, while if in the same example we take
\begin{eqnarray}\label{e5}
c(t)=
  \begin{cases}
     (1+t)^2, ~~ ~0\leqslant t\leqslant 1, &
     \\
     4t, ~~ ~~ ~~~~~~~~ ~~t \geqslant 1,
  \end{cases}
 \end{eqnarray}
then it can be shown that condition ($ii$) of Theorem \ref{th2} is a sufficient condition but not necessary.\hfill$\Box$
\end{r1}
\hspace{0.5in}The following theorem shows that the DAMRL model
preserves the IFRA (increasing failure rate in average) (resp. DFRA (decreasing failure rate in average)) property under certain conditions on $c(t)$. It is useful to remind that a random variable $X$ with failure rate function $r_{X}(\cdot)$ is said to be IFRA  (resp. DFRA) if $\left(\frac{1}{t}\displaystyle\int_{0}^{t}r_{X}(u)du\right)$ is increasing (resp. decreasing) in
$t>0$.
\begin{t1}\label{th3}
If the random variable $X$ is IFRA, then the random variable $X^*$ satisfying (\ref{e1}) is IFRA provided, for all $t\geqslant 0$,
\begin{enumerate}
\item[(i)] $~\frac{m_{X}(t)}{c(t)}$ is increasing in $t$;
\item[(ii)]$~\frac{c'(t)}{c(t)+m_{X}(t)}$ is increasing in $t$.
\end{enumerate}
\end{t1}
$\mathbf{Proof}$: Note, from (\ref{e2}), that
 \begin{eqnarray*}
 \frac{1}{t}\displaystyle\int_{0}^{t}r_{X^*}(u)du&=&\frac{1}{t}\displaystyle\int_{0}^{t}r_{X}(u)\frac{m_{X}(u)}{c(u)+m_{X}(u)}du
 +\frac{1}{t}\displaystyle\int_{0}^{t}\frac{c'(u)}{c(u)+m_{X}(u)}du\\
 &&\\
 &=&P(t),~say.
 \end{eqnarray*}
 Differentiating $P(t)$ with respect to $t$, we get that $P'(t)\geqslant 0$ for all $t\geqslant 0$, if
\begin{enumerate}
\item[(a)] $~~\frac{r_{X}(t) m_{X}(t)}{c(t)+m_{X}(t)}\geqslant\frac{1}{t}\displaystyle\int_{0}^{t}\frac{r_{X}(u)m_{X}(u)}{c(u)+m_{X}(u)}du$,
\item[(b)]$~~\displaystyle\int_{0}^{t}\left(\frac{c'(t)}{c(t)+m_{X}(t)}-\frac{c'(u)}{c(u)+m_{X}(u)}\right)du\geqslant 0$.
\end{enumerate}
 Now, since $X$ is IFRA, it can be shown that (a) holds if
 $$\displaystyle\int_{0}^{t}\left(\frac{m_{X}(t)}{c(t)+m_{X}(t)}-\frac{m_{X}(u)}{c(u)+m_{X}(u)}\right)r_{X}(u)du\geqslant 0,$$
 which is true by ($i$). Further, (b) holds by ($ii$). Hence, the result follows.\hfill$\Box$

\begin{r1}\label{r3}
Note that Example \ref{ex13} can be considered as an application of Theorem \ref{th3}.
Counterexample \ref{co4} can be considered to show that condition
($i$) of Theorem~\ref{th3} is a sufficient condition but not necessary. That condition
($ii$) of Theorem~\ref{th3} cannot be dropped, can be seen by taking
$X$ as a standard exponential random variable, and $c(t)=1/(2+t^2)$, for all $t\geqslant 0$.~\hfill$\Box$
\end{r1}
\hspace{0.5in}The following theorem whose proof is similar to that of Theorem \ref{th3}, shows that the DAMRL model
preserves the
DFRA property under certain conditions on $c(t)$.
\begin{t1}\label{th4}
If the random variable $X$ is DFRA, then the random variable $X^*$ satisfying (\ref{e1}) is DFRA provided, for all $t\geqslant 0$,
\begin{enumerate}
\item[(i)] $~\frac{m_{X}(t)}{c(t)}$ is decreasing in $t$;
\item[(ii)]$~\frac{c'(t)}{c(t)+m_{X}(t)}$ is decreasing in $t$.\hfill$\Box$
\end{enumerate}
\end{t1}
\begin{r1}
An application of the above theorem can be considered by taking $X$ as standard exponential and $c(t)$ as given in (\ref{e5}).\hfill$\Box$
\end{r1}
\begin{r1}\label{r4}
Considering $X$ as in Example \ref{ex1} and $c(t)=\frac{1+t}{2+t}$, $t\geqslant 0$, one can see that condition ($i$) of Theorem~\ref{th4}
is a sufficient condition but not necessary. That condition
($ii$) of Theorem \ref{th4} is a sufficient condition can be seen by considering
\begin{eqnarray*}
c(t)=
  \begin{cases}
     t, ~~~~ ~~~~~~~~~~~0\leqslant t\leqslant 0.5, &
     \\
     1/4+t^2,~~~~0.5\leqslant t\leqslant 1,&
     \\
     2t-3/4, ~~ ~~ ~~~~~~~t \geqslant 1.
  \end{cases}
\end{eqnarray*}
and $X$ as in Example \ref{ex1}.\hfill$\Box$
\end{r1}
\hspace{0.5in}The following theorem shows that the DAMRL model
preserves the NBU (new better than used) (resp. NWU (new worse than used)) property under certain conditions on $c(t)$. It is useful to remind that a random variable $X$ is said to have NBU (resp. NWU) property if $~\bar{F}_{X}(t+x)\leqslant (resp. \geqslant)~\bar{F}_{X}(t)\bar{F}_{X}(x)$,
for all $x,~t\geqslant 0$.

\begin{t1}\label{th5}
If the random variable $X$ is NBU, then the random variable $X^*$ satisfying (\ref{e1}) is NBU, provided, for all $t\geqslant 0$,
\begin{enumerate}
\item[(i)] $~1+\frac{c(t)}{m_{X}(t)}$ is logconvex;
\item[(ii)]$~\frac{c(t)}{m_{X}(t)\left(c(t)+m_{X}(t)\right)}$ is decreasing in $t$.
\end{enumerate}
\end{t1}
$\mathbf{Proof}$: It can be shown that $X^*$ is NBU if and only if, for all $x,~t\geqslant 0$,
$$\displaystyle\int_{x}^{t+x}r_{X}(u)du-\displaystyle\int_{0}^{t}r_{X}(u)du+\displaystyle\int_{x}^{t+x}\frac{c'(u)-r_{X}(u)c(u)}{c(u)+m_{X}(u)}du-
\displaystyle\int_{0}^{t}\frac{c'(u)-r_{X}(u)c(u)}{c(u)+m_{X}(u)}du\geqslant 0.$$
\\
Since $X$ is NBU, it is sufficient to show that,
for all $x,~t\geqslant 0$,
$$\displaystyle\int_{x}^{t+x}\frac{c'(u)-r_{X}(u)c(u)}{c(u)+m_{X}(u)}du-
\displaystyle\int_{0}^{t}\frac{c'(u)-r_{X}(u)c(u)}{c(u)+m_{X}(u)}du\geqslant 0.$$
\\
This holds if
$$\frac{d}{dt}\left[\ln\left(1+\frac{c(t)}{m_{X}(t)}\right)\right]-\frac{c(t)}{m_{X}(t)\left[c(t)+m_{X}(t)\right]}
\mathrm{~~~is~increasing~in~{\it t}}.$$
This is true if the hypotheses are true.\hfill$\Box$
\begin{r1}\label{r5}
Example \ref{ex13} can be taken as an application of Theorem \ref{th5}. By taking $X$ as standard exponential, and $c(t)=\frac{1}{2+t^2}$, $t\geqslant 0$, one can show that condition ($i$) of Theorem \ref{th5} cannot be dropped. That condition ($ii$) of Theorem \ref{th5}  is a sufficient condition but not necessary can be seen by taking $X$ a random variable having mean residual life $m_{X}(t)=\frac{1}{2+t}$, for all $t \geqslant 0$, and $c(t)=\frac{1}{2+t}$, for all $t \geqslant 0$.\hfill$\Box$
\end{r1}
\hspace{0.5in}The following theorem whose proof is similar to that of Theorem \ref{th5} shows that the DAMRL model
preserves the
NWU property under certain conditions on $c(t)$.
\begin{t1}\label{th6}
If the random variable $X$ is NWU, then the random variable $X^*$ satisfying (\ref{e1}) is NWU, provided, for all $t\geqslant 0$,
\begin{enumerate}
\item[(i)] $~1+\frac{c(t)}{m_{X}(t)}$ is logconcave;
\item[(ii)]$~\frac{c(t)}{m_{X}(t)\left(c(t)+m_{X}(t)\right)}$ is increasing in $t$.\hfill$\Box$
\end{enumerate}
\end{t1}
\begin{r1}
Since a concave function is logconcave, condition ($i$) of Theorem \ref{th6} can be replaced by
`$\frac{c(t)}{m_{X}(t)}$ is concave'.
\end{r1}
\begin{r1}\label{r8}
To present an application of Theorem \ref{th6}, one may take $X$ as a standard exponential random variable, and $c(t)=t/(1+t),~t\geqslant 0$.\hfill$\Box$
\end{r1}
\hspace{0.5in}The following counterexample shows that condition
($i$) of Theorem \ref{th6} is a sufficient condition but not necessary.
\begin{co}\label{co8}
Let $X$ be a random variable with survival function
\begin{eqnarray*}
\bar{F}_{X}(t)=
  \begin{cases}
     \frac{1+t^{2}}{(1+t)^{3}}\exp\left[t-t^{2}/2\right]\nonumber,~~~~  0\leqslant t \leqslant 1, &
     \\
     \frac{1}{4}\exp\left[1-t/2\right]\nonumber,   ~~ ~~ ~~~~ ~~~~ ~~~t \geqslant 1.
  \end{cases}
 \end{eqnarray*}
It is shown in Nanda, Das and Balakrishnan \cite{ndb} that $X$ is NWU. Take
\begin{eqnarray*}
c(t)=
  \begin{cases}
     t^2, ~~~~ ~~~~0\leqslant t\leqslant 1, &
     \\
     2t-1, ~~ ~~ ~~~~t \geqslant 1.
  \end{cases}
\end{eqnarray*}
Clearly, $c(t)$ satisfies all the conditions of Lemma \ref{l1}. Now, we see that, for $0\leqslant t\leqslant 1$,
$\ln\left(1+\frac{c(t)}{m_{X}(t)}\right)$
is not concave. Again,
\begin{eqnarray*}
&&\frac{c(t)}{m_{X}(t)\left(c(t)+m_{X}(t)\right)}\\&=&
\begin{cases}
\frac{16t^2(1+t^2)^2\exp(2t-t^2)}{\left(4(1+t)\exp[t-t^{2}/2]+(1+t)^{3}\exp(1/2)\right)
\left(4t^2(1+t^2)\exp[t-t^{2}/2]+4(1+t)\exp[t-t^{2}/2]+(1+t)^{3}\exp(1/2)\right)},~~~~0\leqslant t\leqslant 1, &
\\
\frac{2t-1}{2(2t+1)},~~~~~~~~~~~ ~~~~~~~~~~~~~~~~~~~~~~~~~~ ~~~~~~~~~~~~~~~~~~~~~~~~~~~~~~~~~~~~~~~~~~~~~~~~~~~~~t\geqslant 1
\end{cases}
\end{eqnarray*}
is increasing in $t$. Now, a tedious algebra shows that $X^*$ is NWU. Hence, the condition
($i$) of Theorem \ref{th6} is a sufficient condition but not necessary. \hfill$\Box$
\end{co}
\begin{r1}\label{r6}
Taking $X$ to be a random variable having mean residual life $m_{X}(t)=1+t$, $t\geqslant 0$, and $c(t)=t$, for all $t\geqslant 0$,
it can be shown that the condition ($ii$) of Theorem \ref{th6} is a sufficient condition but not necessary.~\hfill$\Box$
\end{r1}
\hspace{0.5in}The following theorem shows that the model in (\ref{e1}) preserves the NBUFR (new better than used in failure rate) property. It is useful to remind that a random variable $X$ is said to be NBUFR (resp. NWUFR) if $r_{X}(t)\geqslant (resp. \leqslant)~r_{X}(0)$, for all $t>0$. The proof is omitted.
\begin{t1}\label{th7}
If the random variable $X$ is NBUFR, then the random variable
$X^{*}$ satisfying (\ref{e1}) is NBUFR, provided, for all $t\geqslant0$,
\begin{enumerate}
\item[(i)]$~\frac{c(t)}{m_{X}(t)}\leqslant \frac{c(0)}{E(X)}$;
\item[(ii)]$~\frac{c'(t)}{c(t)+m_{X}(t)}\geqslant \frac{c'(0)}{c(0)+E(X)}$. \hfill$\Box$
\end{enumerate}
\end{t1}
\begin{r1}\label{r10}
Let $X$ be a random variable having mean residual life $m_{X}(t)=\frac{1}{2+t}$, $t\geqslant 0$, and $c(t)=\frac{1}{3+t}$, $t\geqslant 0$. Then one can show that condition ($i$) of Theorem \ref{th7} is a sufficient condition but not necessary.\hfill$\Box$
\end{r1}
\begin{r1}\label{r7}
Taking $X$ to be a standard exponential random variable and $c(t)=\frac{1}{2+t^2}$, $t\geqslant 0$, it can be shown that condition ($ii$) of Theorem \ref{th7} cannot be dropped. \hfill$\Box$
\end{r1}
\begin{c1}\label{c7}
If the random variable $X$ is NBUFR, then the random variable $X^*$ satisfying (\ref{e1}) is NBUFR provided, for all $t\geqslant 0$,
\begin{enumerate}
\item[(i)] $~\frac{c(t)}{m_{X}(t)}$ is decreasing in $t$;
\item[(ii)]$~\frac{c'(t)}{c(t)+m_{X}(t)}$ is increasing in $t$. \hfill$\Box$
\end{enumerate}
\end{c1}
\begin{r1}
Example \ref{ex13} can be considered as an application of Theorem \ref{th7} and also of Corollary \ref{c7}. \hfill$\Box$
\end{r1}
\hspace{0.5in}The following theorem shows that the model in (\ref{e1}) preserves the NWUFR property.
\begin{t1}\label{th8}
If the random variable $X$ is NWUFR, then the random variable
$X^{*}$ satisfying (\ref{e1}) is NWUFR, provided, for all $t\geqslant0$,
\begin{enumerate}
\item[(i)]$~\frac{c(t)}{m_{X}(t)}\geqslant \frac{c(0)}{E(X)}$;
\item[(ii)]$~\frac{c'(t)}{c(t)+m_{X}(t)}\leqslant \frac{c'(0)}{c(0)+E(X)}$.\hfill$\Box$
\end{enumerate}
\end{t1}
\begin{ex}\label{ex4}
As an application of Theorem \ref{th8}, one may take $X$ to be a random variable as defined in Example \ref{ex1}, and $c(t)$ as defined in (\ref{e5}). \hfill$\Box$
\end{ex}
\begin{co}\label{co11}
Let $X$ be a random variable as defined in Example \ref{ex1}. Taking $c(t)=\frac{1+t}{2+t}$, $t\geqslant 0$, one can show that condition ($i$) of Theorem \ref{th8} is a sufficient condition but not necessary, where as by taking
\begin{eqnarray}\label{e7}
c(t)=
  \begin{cases}
     t^2,~~~~ ~~ ~~~0\leqslant t\leqslant 1, &
     \\
     2t-1,   ~~~~~~~ ~~t \geqslant 1,
  \end{cases}
 \end{eqnarray}
one can show that condition ($ii$) of Theorem \ref{th8} cannot be dropped.~\hfill$\Box$
\end{co}
\begin{c1}\label{c6}
If the random variable $X$ is NWUFR, then the random variable $X^*$ satisfying (\ref{e1}) is NWUFR provided, for all $t\geqslant 0$,
\begin{enumerate}
\item[(i)] $~\frac{c(t)}{m_{X}(t)}$ is increasing in $t$;
\item[(ii)]$~\frac{c'(t)}{c(t)+m_{X}(t)}$ is decreasing in $t$. \hfill$\Box$
\end{enumerate}
\end{c1}
\begin{r1}
Example \ref{ex4} can be considered as an application of Corollary \ref{c6}.\hfill$\Box$
\end{r1}
\hspace{0.5in}The following theorem shows that the model in
(\ref{e1}) preserves the NBAFR (new better than used in failure rate average) property. It is useful to remind that a random variable $X$ is said to be NBAFR (resp. NWAFR) if $\int_{0}^{t}r_{X}(u)du\geqslant (resp. \leqslant)~t r_{X}(0)$, for all $t>0$.
\begin{t1}\label{th9}
If the random variable $X$ is NBAFR, then the random variable
$X^{*}$ satisfying (\ref{e1}) is NBAFR, provided, for all $t\geqslant0$,
\begin{enumerate}
\item[(i)]$~\frac{c(t)}{m_{X}(t)}\leqslant \frac{c(0)}{E(X)}$;
\item[(ii)]$~\frac{c'(t)}{c(t)+m_{X}(t)}\geqslant \frac{c'(0)}{c(0)+E(X)}$.
\end{enumerate}
\end{t1}
$\mathbf{Proof}$:~$X^*$ is NBAFR if, for all $t\geqslant0$,
\begin{eqnarray}\label{e3}
\frac{1}{t}\displaystyle\int_{0}^{t}\frac{r_{X}(u)m_{X}(u)}{c(u)+m_{X}(u)}du+\frac{1}{t}\displaystyle\int_{0}^{t}\frac{c'(u)}{c(u)+m_{X}(u)}du
\geqslant
\frac{r_{X}(0)E(X)}{c(0)+E(X)}+\frac{c'(0)}{c(0)+E(X)}.
\end{eqnarray}
Since $X$ is NBAFR, by hypothesis ($ii$) we have that (\ref{e3}) holds if
$$\frac{1}{t}\displaystyle\int_{0}^{t}r_{X}(u)\left(\frac{m_{X}(u)}{c(u)+m_{X}(u)}-\frac{E(X)}{c(0)+E(X)}\right)du\geqslant 0,$$
which holds by ($i$). Hence the result follows.\hfill$\Box$
\begin{r1}\label{r12}
Let $X$ be a random variable having mean residual life $m_{X}(t)=\frac{1}{2+t}$, $t\geqslant 0$ and $c(t)=\frac{1}{3+t}$, $t\geqslant 0$. Then one can verify that $X^*$ is NBAFR. Hence, condition ($i$) of Theorem~\ref{th9} is a sufficient condition but not necessary. That condition ($ii$) of Theorem \ref{th9} cannot be dropped, can be seen by considering $X$ as a standard exponential random variable and $c(t)=\frac{1}{2+t^2}$, $t\geqslant 0$.\hfill$\Box$
\end{r1}
\begin{c1}\label{c8}
If the random variable $X$ is NBAFR, then the random variable $X^*$ satisfying (\ref{e1}) is NBAFR provided, for all $t\geqslant 0$,
\begin{enumerate}
\item[(i)] $~\frac{c(t)}{m_{X}(t)}$ is decreasing in $t$;
\item[(ii)]$~\frac{c'(t)}{c(t)+m_{X}(t)}$ is increasing in $t$. \hfill$\Box$
\end{enumerate}
\end{c1}
\begin{r1}
Example \ref{ex13} can be considered as an application of Theorem \ref{th9} and also of Corollary \ref{c8}. \hfill$\Box$
\end{r1}
\hspace{0.5in}The following theorem whose proof is similar to that of Theorem \ref{th9} shows that the model in
(\ref{e1}) preserves the NWAFR property.
\begin{t1}\label{th10}
If the random variable $X$ is NWAFR, then the random variable
$X^{*}$ satisfying (\ref{e1}) is NWAFR, provided, for all $t\geqslant0$,
\begin{enumerate}
\item[(i)]$~\frac{c(t)}{m_{X}(t)}\geqslant \frac{c(0)}{E(X)}$;
\item[(ii)]$~\frac{c'(t)}{c(t)+m_{X}(t)}\leqslant \frac{c'(0)}{c(0)+E(X)}$.\hfill$\Box$
\end{enumerate}
\end{t1}
\begin{r1}\label{r9}
Example \ref{ex1}, with $c(t)$ as given in (\ref{e5}), can be taken as an application of Theorem \ref{th10}.\hfill$\Box$
\end{r1}
\begin{r1}\label{r11}
Example \ref{ex1} with $c(t)=\frac{1+t}{2+t}$, $t\geqslant 0$, shows that condition ($i$) of Theorem \ref{th10} is a sufficient condition but not necessary. That the condition ($ii$) of Theorem \ref{th10} cannot be dropped, can be seen by taking Example \ref{ex1}, with $c(t)$ as defined in (\ref{e7}).\hfill$\Box$
\end{r1}
\begin{c1}\label{c9}
If the random variable $X$ is NWAFR, then the random variable $X^*$ satisfying (\ref{e1}) is NWAFR provided, for all $t\geqslant 0$,
\begin{enumerate}
\item[(i)] $~\frac{c(t)}{m_{X}(t)}$ is increasing in $t$;
\item[(ii)]$~\frac{c'(t)}{c(t)+m_{X}(t)}$ is decreasing in $t$. \hfill$\Box$
\end{enumerate}
\end{c1}
\begin{ex}
Let $X$ be a random variable as defined in Example \ref{ex1}. Clearly $X$ is NWAFR. Take $c(t)$ as defined in (\ref{e5}). In Example \ref{ex4}, we see that both the conditions of Theorem~\ref{th10} are satisfied. Hence, by Corollary \ref{c9}, $X^*$ is NWAFR.\hfill$\Box$
\end{ex}
\section{Conclusion}
In this manuscript, we study the properties of the model $m^{*}(t)=c(t)+m(t)$, which may be considered as a time-dependent additive MRL model. We have given conditions under which $m^{*}(\cdot)$ can be considered to be an MRL function of some random variable. We have also studied the conditions under which the variable $X^*$ (having MRL function $m^*$) belongs to the ageing class $\textrm{C}$ when $X$ belongs to the class $\textrm{C}$. To be specific, let $x\in$ IFR such that $x$ satisfies some property $P$. This means that the set of random variables which belong to IFR class satisfying the property $P$, defines a subclass of IFR class, call it $\textrm{C}$. Then the results studied in this manuscript are the closure properties of $\textrm{C}$ ($\subset$ IFR). With the help of counterexamples we have shown that the ageing classes like IFR etc. are not closed under the model discussed whereas a subclass of each of these well known classes are closed.
\section*{Acknowledgements}
This is a part of the PhD thesis of Suchismta Das submitted at IISER Kolkata in 2014. The financial support from NBHM, Department of Atomic Energy, Government of India (Grant no. 2/48(25)/2014/NBHM(R.P.)/R \& D II/1393 dated February 3, 2015) is acknowledged by Asok K. Nanda.

\end{document}